\documentclass{article}
\usepackage{latexcad}
\usepackage{amsfonts,graphicx}
\begin{document}
\def \version {22 March 2010} 
\def \komm {\LARGE\sf }
\def \kkk {{\LARGE\sf --!!!-- }}
\def \sn {{\Bbb N}}
\def \vp {\varphi}
\def \msk {\medskip}
\def \ssk {\smallskip}
\newtheorem{theo}{Theorem}
\newtheorem{lem}{Lemma}
\newtheorem{obs}[lem]{Observation}
\newtheorem{con}[theo]{Conjecture}
\newtheorem{prop}[lem]{Proposition}
\newtheorem{prob}{Problem}
\title{ 
 List colorings of $K_5$-minor-free graphs with
 special list assignments}
\author{
 Daniel W.
  Cranston\thanks{Department of Mathematics and Applied Mathematics, Virginia Commonwealth University, Richmond, Virginia, USA},
  Anja Pruchnewski\thanks{Faculty of Mathematics and Natural Sciences, Ilmenau University of Technology, Ilmenau, Germany},
 Zsolt Tuza\thanks{Computer and Automation Institute,
   Hungarian Academy of Sciences, Budapest;
 and Department of Computer Science and Systems Technology,
  University of Pannonia, Veszpr\'em, Hungary.
  Research supported in part by the Hungarian Scientific Research
  Fund, OTKA grant 81493.},
   Margit Voigt\thanks{Faculty of Information Technology and Mathematics, University of Applied Sciences, Dresden, Germany}}
\date{\small \version }
\maketitle

\begin{abstract}

A {\it list assignment} $L$ of a graph $G$ is a function that
assigns a set (list) $L(v)$ of colors to every vertex $v$ of $G$.
Graph $G$ is called {\it $L$-list colorable} if it admits a vertex
coloring $\varphi$ such that $\varphi(v)\in L(v)$ for all $v\in
V(G)$ and $\varphi(v)\not=\varphi(w)$ for all $vw\in E(G)$.

The following question was raised by Bruce Richter. Let $G$ be a
planar, 3-connected graph that is not a complete graph. Denoting
by $d(v)$ the degree of vertex $v$, is $G$ $L$-list colorable for
every list assignment $L$ with $|L(v)|=\min \{d(v), 6\}$ for all
$v\in V(G)$?

More generally, we ask for which pairs $(r,k)$ the following
question has an affirmative answer. Let $r$ and $k$ be integers
and let $G$ be a $K_5$-minor-free $r$-connected graph that is not
a Gallai tree (i.e., at least one block of $G$ is neither a
complete graph nor an odd cycle). Is $G$ $L$-list colorable for
every list assignment $L$ with $|L(v)|=\min\{d(v),k\}$ for all
$v\in V(G)$?

We investigate this question by considering the components of
$G[S_k]$, where $S_k:=\{v\in V(G)~ |~ d(v)<k\}$ is the set of
vertices with small degree in $G$. We are especially interested in
the minimum distance $d(S_k)$ in $G$ between the components of
$G[S_k]$.

\end{abstract}

\vfill

\section{Introduction}

In the seventies of the last century, the concept of list
colorings  was introduced independently by Erd\H{o}s, Rubin and
Taylor \cite{ErdosRT79} and by Vizing \cite{Vizing76}. Since then
this topic has been studied extensively by many authors, including
\cite{Bo}-\cite{Voig93}.  In particular, list colorings of planar
graphs have received and continue to receive enormous amounts of
attention;  see, e.g., the surveys \cite{Tuza97,KratoTV99}.

Let $G=(V,E)$ be a graph, let $f:V\longrightarrow \sn$ be a
function, and let $k\geq 0$ be an integer. A {\it list assignment}
$L$ of $G$ is a function that assigns to every vertex $v$ of $G$ a
set (list) $L(v)$ of colors (usually each color is a positive
integer). We say that $L$ is an {\it $f$-assignment} or a
{\it $k$-assignment} if $|L(v)|=f(v)$ for all $v\in V$ or
$|L(v)|=k$ for all $v\in V$, respectively. A {\it coloring} of $G$
is a function $\vp$ that assigns a color to each vertex of $G$ so
that $\vp(v)\not=\vp(w)$ whenever $vw\in E$. An {\it $L$-coloring}
of $G$ is a coloring $\vp$ of $G$ such that $\vp(v)\in L(v)$ for
all $v\in V$. If $G$ admits an $L$-coloring, then $G$ is
{\it $L$-colorable}.
When $L(v)=\{1\dots,k\}$ the corresponding terms
become {\it $k$-coloring} and {\it $k$-colorable}, respectively. The
graph $G$ is said to be {\it $f$-list colorable} if $G$ is
$L$-colorable for every $f$-assignment $L$ of $G$. When $f(v)=k$
for all $v\in V$, the corresponding term becomes {\it $k$-list
colorable}.

Erd\H{o}s, Rubin and Taylor \cite{ErdosRT79} asked, among
other problems, the following two questions.  Are there planar
graphs that are not 4-list colorable? Is every planar graph 5-list
colorable? Both questions were answered in 1993.  In
\cite{Voig93} Voigt gave the first example of a non 4-list
colorable planar graph.  In \cite{Thom94}  Thomassen
answered the second question with a beautiful proof of the
following result.

\begin{theo}[\cite{Thom94}] \label{five}
Every planar graph is 5-list colorable.
\end{theo}

 \v{S}krekovski \cite{Skre98} extended this result to
$K_5$-minor-free graphs.

\begin{theo}[\cite{Skre98}] \label{Skre}
Every\/ $K_5$-minor-free graph is 5-list colorable.
\end{theo}

In 2008 Hutchinson \cite{Hutch} published results on list
colorings of subclasses of planar graphs, where,
for every $v\in V(G)$, the function
$f(v)$ is the minimum of the vertex degree $d(v)$ and a given
integer.

\begin{theo}[\cite{Hutch}]
If\/ $G$ is a 2-connected outerplanar bipartite graph and\/
$f(v)=\min\{d(v),4\}$ for all\/ $v \in V$, then\/ $G$ is\/
$f$-list colorable.
\end{theo}

\begin{theo}[\cite{Hutch}]
If\/ $G$ is a 2-connected outerplanar near-triangulation and\/
$f(v)=\min\{d(v),5\}$ for all\/ $v\in V$, then\/ $G$ is\/ $f$-list
colorable except when the graph is\/ $K_3$ with identical 2-lists.
\end{theo}

In the same paper Hutchinson mentioned the following problem posed
by Bruce Richter.

\begin{prob}[\cite{Hutch}] \label{p1}
Let\/ $G$ be a planar, 3-connected graph that is not a complete
graph and let\/ $f(v)= \min \{d(v), 6\}$ for all\/ $v\in V$. Is\/
$G$ $f$-list colorable?
\end{prob}

In this paper we give partial results concerning the above problem.  Here we study the class of $K_5$-minor-free graphs, which contains the class of planar graphs as a subset. We also investigate the analogous question for non 3-connected $K_5$-minor-free graphs.
In that case both the complete graphs and the so-called Gallai trees play a special role.  A {\it Gallai tree} is a graph $G$ such that every block of $G$ is either a complete graph or an odd cycle. Let $\kappa(G)$ denote the {\it connectivity} of $G$, that is, the cardinality of a smallest vertex cut set of $G$.

\begin{prob} \label{p2}
Let\/ $r\geq 1$ and\/ $k\geq 5$ be integers. Let\/ $G$ be a\/
$K_5$-minor-free graph with $\kappa(G)=r$, such that\/ $G$ is not
a Gallai tree.
Is\/ $G$ $f$-list colorable when\/ $f(v)=\min\{d(v),k\}$?
\end{prob}

An important tool for our investigations is the following theorem.

\begin{theo}[{\cite{Bo,ErdosRT79,KostSW96}}] \label{gallai}
Let\/ $G$ be a connected graph and let\/ $L$ be a list assignment
with\/ $|L(v)|\geq d(v)$ for all\/ $v \in V$. If\/ $G$ has no\/
$L$-coloring, then the following three conditions hold.
\begin{enumerate}
\item[(a)] $|L(v)|=d(v)$ for every vertex\/ $v\in V(G)$.

\item[(b)] G is a Gallai tree.

\item[(c)] Let\/ $\cal B$ be the set of blocks of\/ $G$ and\/
 ${\cal B}(v)\subseteq{\cal B}$ the set of blocks of\/ $G$
  containing a specified vertex\/ $v$. There exist
 color sets\/ $L(B)$ for all\/ $B \in {\cal B}$ such that\/
  $L(B_1) \cap L(B_2) =\emptyset$ whenever\/ $B_1$ and\/
   $B_2$ have a common vertex and for all\/ $v\in V(G)$ we have\/
    $L(v) =\bigcup\limits_{{\cal B}(v)} L(B)$.
\end{enumerate}
\end{theo}

In this paper we investigate Problem \ref{p2}, considering subsets
of planar graphs that fulfill special requirements. Let
$$S_k:=\{v\in V(G)~ |~ d(v)<k\}~~ {\rm and}~~ B_k:=\{v\in V(G)~ |~
d(v)\geq k\}$$ be the sets of vertices with small degree in $G$
and with big degree in $G$, respectively. The smallest distance
between components of $G[S_k]$ in $G$ is denoted by $d(S_k)$,
where $G[S_k]$ is the subgraph of $G$ induced by $S_k$. If
$G[S_k]$ has at most one component, then let $d(S_k)=0$.
We may always assume that $G$ is connected since otherwise we can consider seperately each component of $G$.
We will answer the question of Problem \ref{p2} for many cases. Our results for $K_5$-minor-free graphs are summarized in the following tables.

\medskip

\renewcommand{\arraystretch}{1.3}

\begin{center}
\begin{minipage}{6cm}
{\quad\bf Table 1:}
$\kappa(G)\in\{1,2\}$ \\

\begin{tabular}{|r|c|c|c|c|}
 \cline{1-5}
{\small $_k \backslash ^{d(S_k)}$} & 2 & 3 & 4 & $\geq 5$ \\
\cline{1-5} 5 & -- & -- & -- & ? \\ \cline{1-5} 6 & -- & --/? & ?
& +
\\ \cline{1-5} 7 & -- & --/? & ? & + \\ \cline{1-5} $\geq 8$ & -- & +
& + & + \\ \cline{1-5}
\end{tabular}
\end{minipage}
\begin{minipage}{5cm}
{ ~ \,\bf Table 2:} $\kappa(G)\in\{3,4\}$\\

\begin{tabular}{|r|c|c|c|c|}
 \cline{1-5}
{\small $_k \backslash ^{d(S_k)}$} & 2 & 3 & 4 & $\geq 5$ \\
\cline{1-5} 5 & -- & -- & -- & ? \\ \cline{1-5} 6 & {--} & { ?} &
{ ?} & + \\ \cline{1-5} 7 &  -- & + & + & + \\ \cline{1-5} $\geq
8$ & -- & + & + & + \\ \cline{1-5}
\end{tabular}
\end{minipage}\end{center}

\medskip

In Section 2 we give some results for $k\geq 6$ that can be obtained by simple observations including the solution for $\kappa(G)\geq 5$. Section 3 contains our main results for connected graphs, whereas Section 4 deals with graphs with $\kappa(G)\in\{3,4\}$. In Section 5 we consider the case $k=5$ and in Section 6 we mention open problems.

Note that the original problem asked for planar graphs with $\kappa(G)\in\{3,4\}$. In that case the answer for $d(S_k)=2$ and $k\geq 6$ is still unknown. All other entries of the above tables are valid also for planar graphs.

\section{Observations}

In this section we collect some immediate results. Let $k=6$, that
is, let $$f(v)=\min\{d(v),6\}~~ {\rm for~ all}~~ v\in V.$$

\begin{obs}
If\/ $G$ is a\/ $K_5$-minor-free graph and\/ $d(v)\geq 5$ for
all\/ $v\in V$, then\/ $G$ is\/ $f$-list colorable.
\end{obs}

In this case we have $|L(v)|\geq 5$ for all $v\in V$ and we are
done since every $K_5$-minor-free graph is 5-list colorable
\cite{Skre98}. If $G$ is 5-connected, then the degree of each
vertex is at least 5. So our next observation follows immediately.

\begin{obs}
If\/ $G$ is a 5-connected\/ $K_5$-minor-free graph, then\/ $G$
is\/ $f$-list colorable.
\end{obs}

The next observation follows from Theorem \ref{gallai}.
\begin{obs}
If\/ $G$ is a\/ $K_5$-minor-free graph that is not a Gallai tree
and\/ $d(v)\leq 6$ for all\/ $v\in V$, then\/ $G$ is\/ $f$-list
colorable.
\end{obs}

\begin{obs}
Let\/ $G$ be a\/ $K_5$-minor-free graph. If the vertices of degree
at most 5 have pairwise distance at least 3 in\/ $G$, then\/ $G$
is\/ $f$-list colorable. \label{dist3}
\end{obs}

To prove Observation~\ref{dist3}, we color each vertex $v\in S_6$
with an arbitrary color from its list.  We delete the color used
on each $v$ from the lists of the neighbors of $v$ and remove the
colored vertices from $G$, obtaining $G^*$. Because of the
hypothesis, every remaining vertex has at most one neighbor $v\in
S_6$ in $G$. Thus we have $|L^*(v)|\geq 5$ for the reduced lists
of the vertices of $G^*$, so $G^*$ is $L^*$-list colorable by
Theorem \ref{Skre}.

\begin{prop}
If\/ $G$ is a\/ $K_5$-minor-free graph with\/ $d(S_6)\geq 5$,
then\/ $G$ is\/ $f$-list colorable.
\end{prop}
{\bf Proof.} Let $G$ be a $K_5$-minor-free graph with $d(S_6)\geq 5$ and let $L$ be a list assignment with $|L(v)|=\min\{d(v),6\}$ for all $v \in V$. Let $G_1, G_2,\dots$ be the components of $G[S_6]$. For each $G_i$, choose a vertex $w_i \in V\setminus S_6$  that has at least one neighbor $v_i$ in $G_i$.
Color each $w_i$ by  $\varphi_i\in L(w_i)\setminus L(v_i)$.
This is always possible since
$ |L(w_i)| = d(w_i)\geq 6 > d(v_i) = |L(v_i)|$.
Delete $\varphi_i$ from the lists of the neighbors of $w_i$ in $G$, obtaining a list assignment $L'$.
Form subgraph $G'$ from $G$, by removing the vertices of $S_6$ and the vertices $w_i$.
%

Since the vertices $w_i$ have distance at least 3 to each other in $G$,
each vertex $v$ of $G'$ still has at least 5 colors available in $L'$. Thus, we can list color $G'$ from the list assignment $L'$ by Theorem \ref{Skre}. Now consider each $G_i$. 
For every vertex $v\in V(G_i)$, delete the colors from its list used on its colored neighbors. We obtain a new list assignment $L''$. Note that $|L''(v)|\geq d(v)$ for all $v\in V(G_i)$ and $|L''(v_i)|>d(v_i)$. By Theorem \ref{gallai}, each $G_i$ is $L''$-list colorable.
This completes  the $L$-list coloring of $G$. \hfill $\Box$

\medskip

Note that an analogous proof works if $G[S_6]$ has only one component.

Now we consider $f(v)=\min\{d(v),k\}$  for arbitrary $k\geq 3$.
If $G$ is at most 2-connected, then there are graphs that are not $f$-list colorable. The first examples, with $\kappa(G)=1$, were given in \cite{Hutch}. Here we give an example with $\kappa(G)=2$ and minimum degree 3.

\begin{prop}
Let\/ $k\geq 3$ be an integer. There are planar 2-connected
non-complete graphs\/ $G$ with\/ $d(S_k)=2$ and\/ $\delta(G)=3$
and list assignments\/ $L$ with\/ $|L(v)|=\min \{d(v), k\}$ such
that\/ $G$ is not\/ $L$-list colorable.
\end{prop}

{\bf Proof.}
Let $L$ be the following list assignment for the graph $G$ in Figure 1, with $s={k\choose 2}$.
\begin{itemize}
\item $L(x)=L(y)=\{1,\dots,k\}$
\item $L(u_i)=L(v_i)=L_i$, where
 $\{L_1,\dots,L_s\}=\{\{i,j,0\}~ |~ 1\leq i < j\leq k\}$
\end{itemize}

\begin{figure}[h]
\unitlength1.3mm
\begin{center}
\begin{picture}(60,50)
\thinlines
\drawpath{6.0}{8.0}{46.0}{8.0}
\drawpath{6.0}{8.0}{26.0}{14.0}
\drawpath{56.0}{14.0}{56.0}{14.0}
\drawpath{26.0}{14.0}{26.0}{18.0}
\drawpath{26.0}{18.0}{6.0}{8.0}
\drawpath{26.0}{14.0}{46.0}{8.0}
\drawpath{46.0}{8.0}{26.0}{18.0}
\drawpath{6.0}{8.0}{26.0}{36.0}
\drawpath{26.0}{36.0}{26.0}{42.0}
\drawpath{26.0}{42.0}{6.0}{8.0}
\drawpath{26.0}{36.0}{46.0}{8.0}
\drawpath{46.0}{8.0}{26.0}{42.0}
\drawcenteredtext{6.0}{6.0}{$x$}
\thicklines
\drawcenteredtext{46.0}{6.0}{$y$}
\drawcenteredtext{26.0}{12.0}{$u_1$}
\drawcenteredtext{26.0}{20.0}{$v_1$}
\drawcenteredtext{26.0}{32.0}{$u_s$}
\drawcenteredtext{26.0}{44.0}{$v_s$}
\thinlines
\drawdotline{26.0}{30.0}{26.0}{22.0}
\end{picture}
\end{center}\vspace{-1cm}
\caption{$G$ is planar and 2-connected, with $\delta(G)=3$.}
\end{figure}
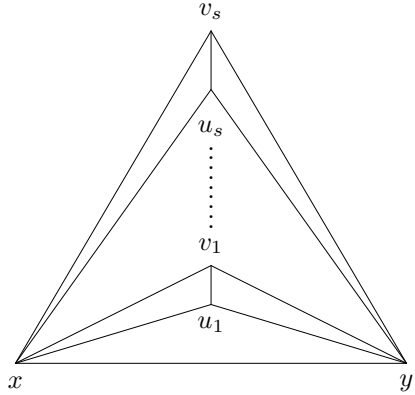

Since the vertex pair $x,y$ has to be colored by one of the pairs $i,j$ with $1\leq i < j\leq k$, it follows that one of the vertex pairs $u_\ell,v_\ell$ with $1\leq \ell \leq s$ is not colorable, since the only remaining available color is $0$.
\hfill $\Box$

\medskip

The above proposition shows that in the original Problem 1  we cannot
replace the assumption of 3-connectivity by 2-connectivity and/or
$\delta(G)=3$.

Furthermore, Hutchinson \cite{Hutch09} pointed out that in Problem 1, even if we
keep 3-connectivity, the requirement of planarity cannot be
dropped; what is more, it cannot be relaxed to the restriction of
being $K_5$-minor-free.

\begin{prop}  \label{c3,d2}
Let\/ $k\geq 3$ be an integer. There are 3-connected,
non-complete,\/ $K_5$-minor-free graphs\/ $G$ with\/ $d(S_k)=2$
and list assignments\/ $L$ with\/ $|L(v)|=\min \{d(v), k\}$ such
that\/ $G$ is not\/ $L$-list colorable.
\end{prop}

{\bf Proof.} Let $s={k\choose 3}$.  Let $G\cong K_{s+3}-E(K_s)$ be
the graph with $V(G)=\{v_1,v_2,\ldots,v_s\}\cup\{x,y,z\}$ such
that $x$, $y$, and $z$ are pairwise adjacent and every $v_i$ is
adjacent to each of $x$, $y$, and $z$. Let $L$ be the following
list assignment.
\begin{itemize}
\item $L(x)=L(y)=L(z)=\{1,\dots,k\}$
\item $L(v_i)=L_i$, where
 $\{L_1,\dots,L_s\}=\{\{h,i,j\}~ |~ 1\leq h < i < j\leq k\}$
\end{itemize}

The three vertices $x$, $y$, $z$ have to be colored by one of the triples $h,i,j$ with $1\leq h < i < j\leq k$.  It follows that one of the vertices $v_\ell$ with $1\leq \ell \leq s$ is not colorable, since the three colors in its list are used on $x$, $y$, and $z$.
\hfill $\Box$

\section{Small connectivity}

\begin{theo}
If\/ $5\leq k\leq 7$ and\/ $f(v)=\min\{d(v),k\}$, then there are
planar non-complete graphs\/ $G$ with\/ $d(S_k)=3$ that are not\/
$f$-list colorable.
\end{theo}

{\bf Proof.} Let $G'$ be a minimal non-4-list colorable
planar graph and $L'$ a 4-assignment such that $G'$ is not
$L'$-list colorable. Thus $\delta(G)\geq 4$. Let $V(G')=\{v_1,\dots,v_n\}$. Take
$n$ copies of $K_3$ where the vertices of the $i$-th copy are
denoted by $x_i,y_i,z_i$. Build $G$ by joining each $v_i$ with
$x_i, y_i$ and $z_i$. Thus $d_G(v_i)\geq 7$.

Let $1,2,3$ be colors that are not contained in the lists $L'(v_i)$ for $i=1,\dots,n$. Let $L(v_i)=L'(v_i)\cup\{1,2,3\}$ and $L(x_i)=L(y_i)=L(z_i)=\{1,2,3\}$. Thus $|L(v_i)|=7$ and $|L(x_i)|=|L(y_i)|=|L(z_i)|=d(x_i)=d(y_i)=d(z_i)=3$.

It is easy to see that $G$ is not $L$-list colorable. Obviously the same construction with shorter lists works if $k=5$ or $k=6$.
\hfill $\Box$

\vspace{1ex}

Note that our example is a graph $G$ with $\kappa(G)=1$. We do not
know an analogous example for $\kappa(G)=2$.  However, the example
shows that, at least for small $k$, the additional assumption
$d(S_k)\geq 3$ is not sufficient if we do not have an assumption on the degree of connectivity.  In contrast, we show next that if
$k\geq 8$, then a list coloring is always possible.

\begin{theo} \label{eight}
Let\/ $k\geq 8$ be an integer and let\/ $f(v)=\min\{d(v),k\}$ for
all\/ $v\in V$. If\/ $G$ is a\/ $K_5$-minor-free graph with\/
$d(S_k)\geq 3$ such that no component of\/ $G$ is a Gallai tree,
then\/ $G$ is\/ $f$-list colorable.
\end{theo}

{\bf Proof.} We prove the theorem for $k=8$. Assume the theorem is
false. Let $G$ be a smallest counterexample and let $L$ be a list
assignment such that $G$ is not $L$-list colorable. Since $G$ is a
smallest counterexample it must be connected.

If all vertices of $G$ have degree at most 8, then $d(v)=|L(v)|$
for all $v\in V$. Theorem \ref{gallai}
shows that such a $G$ is $L$-list colorable if it is not
a Gallai tree. This is a hypothesis of the present theorem, so we
are done. If all vertices of $G$ have degree at least 5, then we
are also done, since every $K_5$-minor-free graph is 5-list
colorable \cite{Skre98}. Thus,  we may assume that both $S_8$
and $B_8$ are nonempty. \msk

Now we would like to apply the following strategy. Let $H$ be a
component of $G[S_8]$. If $H$ is not a Gallai tree, then we remove
it from the graph. If $H$ is a Gallai tree, then we color some
vertices of $H$ and delete the used colors from the lists of all
corresponding neighbors in $G$ (see the case distinction below).
If such a neighbor $w$ belongs to $B_8$, then we shall delete at
most three colors from its list. After that, we remove $H$. If
all components of $G[S_8]$ are removed, then the remaining
subgraph is $K_5$-minor-free and all its lists have cardinality at
least $8-3=5$, since $d(S_8)\geq 3$ in $G$. By Theorem~\ref{Skre}
the remaining graph $G[B_8]$ is list colorable from its reduced
lists.

Now we  re-insert all uncolored vertices and we delete from
their lists all colors used for corresponding neighbors. Let $H'$
denote the subgraph of $H$ induced by the set of all uncolored
vertices of $H$. We shall show that each $H'$ is list colorable
from the reduced lists $L'$. For all uncolored vertices, we know
that $|L'(v)|\geq d_{H'} (v)$. Of course, it suffices to show that
either $G$ is not a smallest counterexample or
$V(H')=\emptyset$.  By Theorem~\ref{gallai}, it also suffices to
show in each case that $H'$ is connected and that it satisfies one
of the following three conditions:

\begin{itemize}
\item[(i)] $H'$ is not a Gallai tree or
\item[(ii)] there is a $v$ in $V(H')$ such that $d_{H'}(v)<|L'(v)|$ or
\item[(iii)] there are adjacent vertices $v$ and $w$ in $H'$ that are not cut vertices and that satisfy $L'(v)\not= L'(w)$.
\end{itemize}
If this strategy works, then the graph $G$ is $L$-list colorable,
which contradicts the assumption that it is a counterexample to the theorem.
\bigskip

If $H$ is not a Gallai tree, then we let $H'=H$, and we are done
by (i). If $H$ is a Gallai tree, then  its blocks are odd
cycles and/or complete graphs on at most four vertices, since $G$
is $K_5$-minor-free. We distinguish the following seven cases.

\msk

{\bf (1) \boldmath $H$ is a complete graph $K_l$, $l\leq 4$.}

\msk

The absence of a $K_5$-minor implies that each vertex in $V(G)$ has at most 3 neighbors
in $H$.  Since $G$ is connected, at least one vertex $v$ of $H$
satisfies the inequality $|L(v)|>d_{H}(v)$. Thus we can
color $H$ properly from its lists. Note that we use at most 3
colors from the list of each neighbor in $G$. Now
$V(H')=\emptyset$, so we are done.

\msk

{\bf (2) \boldmath $H$ is an odd cycle $C_{2l+1}$, $l\geq 2$.}

\msk

Denote the vertices of the cycle by $v_1,\dots, v_{2l+1}$. If
there are adjacent vertices $v_i, v_j$ on the cycle such that
there exists a color $c \in L(v_i)\setminus L(v_j)$, then we use
color $c$ on $v_i$. Now $H'$ is a path, where $v_j$ is an end
vertex, and $|L'(v_j)|>d_{H'} (v_j)$, so we are done by (ii).

Otherwise we have $L(v_1)=\dots=L(v_{2l+1})$. Since $G$ is
connected, $|L(v_i)|\geq 3$ for all $i$. We color the cycle with 3
colors from the lists. Thus $V(H')=\emptyset$.

\msk

{\bf (3) \boldmath $H$ has an end block with a vertex $v$ that is
a cut vertex of~$G$.}

\msk

Let $H_1$ be the end block.  Since $v$ is a cut vertex of $G$ and $|L(u)|=d_H(u)$ for all $u\in (V(H_1)\setminus\{v\})$, we can color all the vertices of $H_1\setminus\{v\}$ before we color $G[B_8]$.  Since $G$ is a smallest counterexample, we can color $G\setminus (H_1\setminus\{v\})$ from its lists.

\msk




{\bf (4) \boldmath $H$ has an end block that is an odd cycle
$C_{2l+1}$, $l\geq 2$.}

\msk

Denote the vertices of the cycle by $v_1,\ldots,v_{2l+1}$, and let
$v_1$ be the cut vertex. Since we are not in case (3) or (iii), we
know that for all $i, j\in \{2,3,\ldots,2l+1\}$ we have
$|L(v_i)|\ge 3$ and $L(v_i) = L(v_j)$.  Thus, we can color $v_2$
through $v_{2l+1}$ with 3 colors, so that $v_2$ and $v_{2l+1}$ get
the same color.  Now $d_{H'} (v_1) < |L'(v_1)|$ and $H'$ is
connected, so we are done by (ii). \msk

{\bf (5) \boldmath $H$ has an end block that is a $K_2$.}
\msk

Denote the vertices of the end block by $v_1$ and $v_2$, where $v_1$ is the cut vertex.
Since we are not in case (3), we have $d_G(v_2)\ge 2$.
Choose two colors, say $a$ and $b$, from $L(v_2)$ and delete them from the lists of all neighbors of $v_2$ belonging to $G[B_8]$. Remove $H$ from $G$. After the coloring of $G[B_8]$, add $H'=H$ and assign $L'(v_2)=\{a,b\}$. Since $d_{H'} (v_2)=1$, we are done by (ii).
\msk

{\bf (6) \boldmath $H$ has an end block that is a $K_3$.}
\msk

Denote the vertices of the end block by $v_1, v_2$ and $v_3$, where $v_1$ is the cut vertex.
Since we are not in case (3), at least one of $v_2$ and $v_3$ has a list of cardinality at least 3, say $|L(v_2)|\geq 3$. Choose 3 colors from this list for $v_2$ and continue as in the previous case.
\ssk

Note that  the same
approach would work for endblocks $K_4$ if we wanted to prove the
theorem only for $k\ge 9$, rather than $k\ge 8$. However, to
prove the theorem for $k=8$, we must consider this final case.

\msk

{\bf (7) \boldmath All end blocks of $H$ are $K_4$.}
\msk

Fix one of the end blocks and denote its vertices by $v_1, v_2, v_3$ and $v_4$, where $v_1$ is the cut vertex.
If there is a color $a\in L(v_i)\setminus L(v_j)$ where $\{i,j\}\subseteq \{2,3,4\}$ then color $v_i$ by $a$ and continue as prescribed in the strategy. It follows that $H'$ is connected and $d_{H'} (v_j) < |L'(v_j)|$, so we are done by (ii). \ssk

Thus $L(v_2)=L(v_3)=L(v_4)$. If $|L(v_i)|=3$ for $i=2,3,4$ then we are done by case (3).
So we know that $|L(v_i)|\geq 4$ for all $i$, and each of $v_2,v_3$ and $v_4$ has at least one neighbor belonging to $B_8$.\ssk

\ssk

{\bf (7a) \boldmath There are $\{i,j\}\in\{2,3,4\}$ and a vertex
$w\in B_8$ such that $w$ is a neighbor of $v_i$, but not a
neighbor of $v_j$. }

\ssk

Remove $H$ from $G$ and join $w$ with all neighbors of $v_j$
belonging to $B_8$. Note that the new graph is still
$K_5$-minor-free. After removing all components of $G[S_8]$, color
the remaining graph from the corresponding lists. This is possible
since the graph is $K_5$-minor-free and all lists have cardinality
at least 5. Remove the additional edges, add $H$ and delete the
colors of the neighbors of the vertices of $H$ from their lists.
Let $a$ be the color used on $w$.  If $a\not\in L(v_i)$, then we
have $d_{H'}(v_i)<|L'(v_i)|$ and we are done by (ii). Otherwise
$a\in L(v_i)=L(v_j)$. But no neighbor of $v_j$ is colored by $a$
since all of them were joined to $w$ by the additional edges.
Since $a\not\in L'(v_i)$ and $a\in L'(v_j)$, we have $L'(v_i)\not=
L'(v_j)$, so we are done by (iii).

\ssk

{\bf (7b) \boldmath Each vertex $z\in B_8$ that is adjacent to at
least one of $v_2,v_3,v_4$ is adjacent to all of them.}

\ssk

Note that $v_2,v_3,v_4$ must have some neighbors in $B_8$ since we
are not in case (3). Let $z\in B_8$ be such a neighbor. Then there
is no path from $v_1$ to $z$ that does not use $v_2,v_3$ or $v_4$,
since otherwise we have a subdivision of $K_5$. This statement
holds also for all further neighbors of $v_2,v_3,v_4$ in $B_8$. So
by removing the edges $v_1v_2, v_1v_3$ and $v_1v_4$, we get two
components $G_1$ containing $v_1$ and $G_2$ containing
$v_2,v_3,v_4$ and
 $z$.

If $G_1$ is a Gallai tree, then we can color it since
$|L(v_1)|>d_{G_1}(v_1)$. Otherwise we can color it, since it
satisfies the hypothesis of the theorem and $G$ is a smallest
counterexample to it. Delete the color of $v_1$ from the lists of
$v_2,v_3$ and $v_4$ and consider $G_2$ with the reduced list
assignment $L'$. We need to prove the following claim to argue
that $G_2$ satisfies the hypothesis of the theorem.

\msk

{\bf Claim:} $G_2$ {\it is not a Gallai tree.}

\ssk

 Assume to the contrary  that $G_2$
is a Gallai tree. Note that $G_2$ has at least one vertex of
degree at least 8, namely $z$. Let $u$ be a vertex of degree at
least 8 that has the largest distance from $v_2$ in $G_2$. This
$u$ belongs to at least 3 blocks of $G_2$, since it has at most 3
neighbors in each block (as $G$ is $K_5$-minor-free). Moreover,
$u$ can have neighbors of degree less than 8 in at most one of
these blocks, since $d(S_8)\geq 3$. Thus there are vertices of
degree at least 8 that have a greater distance from $v_2$ than
$u$, contradicting the assumption for $u$. Hence the claim is
proved. \msk

Thus $G_2$ satisfies the hypothesis of the present theorem and we can color it from the reduced lists, since $G$ is a smallest counterexample. By combining the colorings of $G_1$ and $G_2$, we obtain an $L$-list coloring of $G$, which is a contradiction.

This completes the proof of the theorem.
\hfill $\Box$\msk

\section{$G$ is 3-connected}

\begin{theo} \label{seven}
Let\/ $k\geq 7$ be an integer and let\/ $G$ be a\/
$K_5$-minor-free, non-complete, 3-connected graph. If\/
$d(S_k)\geq 3$, then\/ $G$ is\/ $f$-list colorable when\/
$f(v)=\min \{d(v), k\}$.
\end{theo}

{\bf Proof.}
We use a strategy similar to the proof of Theorem \ref{eight}. Let $G$ be a smallest counterexample to the theorem.
Since $G$ is 3-connected and non-complete, $G$ is not a Gallai tree.
So, if all vertices of $G$ have degree at most 7, then we are done by Theorem~\ref{eight}.  Similarly, if all vertices of $G$ have degree at least 5 then we are done by Theorem \ref{Skre}.

For each component $H$ of $G[S_7]$, we would like to color one or more vertices of $H$ such that for each adjacent vertex in $B_7$ we have to delete at most two colors from its list. If we do this for all components of $G[S_7]$, then we can color $G[B_7]$ from the reduced lists. This is possible because $d(S_7)\geq 3$, which ensures that all reduced lists have cardinality at least 5. Note that we do not need a connectivity assumption for this argument. For a component $H$ of $G[S_7]$, let $H'$ be the subgraph induced by the uncolored vertices of $H$. For every vertex $v\in V(H')$, we obtain a reduced list $L'(v)$ by deleting all colors from $L(v)$ that were used for the already colored neighbors of $v$. Finally, we try to color every $H'$ from the list assignment $L'$.  It will suffice to show, for each $H'$, that $H'$ is connected and satisfies one of conditions (i), (ii) and (iii), from Theorem~\ref{eight}.\msk

If $H$ is not a Gallai tree, then let $H'=H$, and we are done by (i). Thus we can assume that $H$ is a Gallai tree.\msk

{\bf (1) \boldmath $H$ is a single vertex $v$. }
\msk

Color $v$. Since $V(H')=\emptyset$, we are done.\msk

{\bf (2) \boldmath $H$ is $K_2$ or $H$ has an end block $K_3$.} \msk

If $H$ is $K_2$, then denote the vertices of $H$ by $v_1$ and
$v_2$. If $H$ has an end block $K_3$, then denote the vertices of
the end block by $v_1,v_2$ and $v_3$, where $v_3$ is the cut
vertex. If there is an $a\in L(v_i)\setminus L(v_j)$
($\{i,j\}=\{1,2\}$), then color $v_i$ by $a$. In this case,
$|L'(v_j)|>d_{H'}(v_j)$, so we are done by (ii).  So instead, we
may assume $L(v_1)=L(v_2)$.

\ssk

{\bf (2a) \boldmath There is a vertex $z\in B_7$ that is adjacent
to $v_i$ and not adjacent to $v_j$ ($\{i,j\}=\{1,2\}$).}

\ssk

Remove $H$ and join $z$ with all neighbors of $v_j$. The new graph
is still $K_5$-minor-free. Let $H'=H$, with the lists reduced by
the coloring of $G[B_7]$ with the additional edges. Let $a$ be the
color of $z$ in this coloring. If $a\not\in L(v_i)$, then we have
$|L'(v_i)|>d_H(v_i)$, so we are done by (ii). Otherwise we have
$a\not\in L'(v_i)$ and $a\in L'(v_j)$, so we are done by (iii).

\ssk

{\bf (2b) \boldmath  Every $z\in B_7$ that is a neighbor of $v_1$
or $v_2$ is a neighbor of both.}

\ssk

Since $G$ is 3-connected, there have to be at least 3 such
neighbors (including $v_3$ if $H$ is an end block $K_3$). Let
$Z=\{z_1, z_2, z_3,\dots\}$ be the set of all these neighbors.

\ssk

Consider $G'=G[V\setminus\{v_1,v_2\}]$. This subgraph is connected
since $G$ is 3-connected. Thus any vertex pair $z_i,z_j$ is joined
by a path in $G'$. Consider an auxiliary graph $T$ with $V(T)=Z$,
where two vertices $z_i$ and $z_j$  are joined by an edge if and
only if there is path from $z_i$ to $z_j$ in $G'$ not containing
another vertex $z_\ell\in Z$. Note that $T$ has to be a tree since
otherwise there is a subdivision of $K_5$ in $G$.

Let $z_1$ be a leaf of $T$ and $z_2$ its neighbor in $T$. Since
$z_1\in B_7$ or $z_1$ is a cut vertex of $H$, it must have a
neighbor $w\not\in Z\cup\{v_1,v_2\}$.

Now consider a path $P$ from $w$ to $z_3$ in $G$. Since $z_1$ is a
leaf in $T$, every path from $w$ to $z_3$ must contain a vertex
$u\in Z\cup\{v_1,v_2\}$ as an internal vertex. If $P$ contains
$x\in\{v_1,v_2\}$, then, in addition to $z_3$, $P$ must contain
another neighbor of $x$ that is in $Z$; call it $z_i$. Thus every
path from $w$ to $z_3$ contains some vertex $z_i\not=z_3$ as an
internal vertex. Since there have to be at least three internally
disjoint paths from $w$ to $z_3$, we have a path from $z_1$ to
some $z_j\ne z_2$ that has no internal vertices in
$Z\cup\{v_1,v_2\}$.  However, now $T$ is not a tree.  This again
gives us a subdivision of $K_5$, which contradicts the
$K_5$-minor-freeness of $G$.\msk

%

{\bf (3) \boldmath $H$ is an odd cycle $C_{2l+1}$, $l\ge 2$, or $H$ has an end block that is an odd cycle $C_{2l+1}$, $l\ge 2$.} \msk

Denote the vertices of $C_{2l+1}$ by $v_1,v_2, \ldots, v_{2l+1}$, where
the cut vertex (if it exists) is $v_{2l+1}$.
If there is an $a\in L(v_i)\setminus L(v_j)$ ($\{i,j\}\subset\{1,2,\ldots, 2l\}$), then color $v_i$ by $a$. In this case, $|L'(v_j)|>d_{H'}(v_j)$, and we are done by (ii).
So instead, we have $L(v_i)=L(v_j)$ for all $\{i,j\}\subset\{1,2,\ldots,2l\}$. \ssk

If there is a vertex $w\in B_7$ that is adjacent to $v_i$ and not
adjacent to $v_j$ ($\{i,j\}\subset\{1,2,\ldots
2l\}$), then there exist $w'\in B_7$ and $v_{i'}$ and $v_{i'+1}$ ,
which are adjacent on the odd cycle, such that $w'$ is adjacent to
$v_{i'}$ and not adjacent to $v_{i'+1}$. Remove $H$ and join $w'$
with all neighbors of $v_{i'+1}$. The new graph is still
$K_5$-minor free. Let $H'=H$, with the lists reduced by the
coloring of $G[B_7]$ with the additional edges. Let $a$ be the
color of $w'$ in this coloring. If $a\not\in L(v_{i'})$, then
$|L'(v_{i'})|>d_H(v_{i'})$, so we are done by (ii). Otherwise
$a\not\in L'(v_{i'})$ and $a\in L'(v_{i'+1})$, so we are done by
(iii).

Thus we may assume that every vertex  $w\in B_7$ that is adjacent
to  at least one $v_i$, $i\in\{1,\ldots,2l\}$, is adjacent to all
of them. Since $G$ is 3-connected, there must be at least two such
neighbors that lie in the same component of $G[B_7]$. To build a
$K_5$-minor, we take these two neighbors together with $v_1,
v_2$ and $v_3$. \msk

{\bf (4) \boldmath $H$ is $K_3$ or $K_4$ or has an end block that is $K_4$.} \msk

If $H$ is $K_3$, then denote the vertices of $H$ by $v_1,v_2$ and $v_3$.
Otherwise, denote the vertices by $v_1,v_2,v_3$ and $v_4$, where the cut vertex (if it exists) is $v_4$.
If there is an $a\in L(v_i)\setminus L(v_j)$ ($\{i,j\}\subset\{1,2,3\}$), then color $v_i$ by $a$. In this case, $|L'(v_j)|>d_{H'}(v_j)$, and we are done by (ii).
So instead, we have $L(v_1)=L(v_2)=L(v_3)$. \ssk

If there is a vertex $w\in B_7$ that is adjacent to $v_i$ and not adjacent to $v_j$ ($\{i,j\}=\{1,2,3\}$), then remove $H$ and join $w$  with all neighbors of $v_j$. The new graph is still $K_5$-minor free.
Let $H'=H$, with the lists reduced by the coloring of $G[B_7]$ with the additional edges. Let $a$ be the color of $w$ in this coloring. If $a\not\in L(v_i)$, then $|L'(v_i)|>d_H(v_i)$, so we are done by (ii). Otherwise $a\not\in L'(v_i)$ and $a\in L'(v_j)$, so we are done by (iii).

Thus we may assume that every vertex  $w\in B_7$ that is adjacent to  $v_1, v_2$ or $v_3$ is adjacent to all of them. Since $G$ is 3-connected, there must be at least 3 such neighbors that all lie in the same component of $G\setminus \{v_1,v_2,v_3\}$. To build a $K_5$-minor, we take two of these neighbors and $v_1, v_2$ and $v_3$.   \msk

%
%
%

%

{\bf (5) \boldmath All end blocks of $H$ are $K_2$.} \msk

Denote the vertices of an end block by $v_1$ and $v_2$, where $v_1$ is the cut vertex. 
Since $G$ is 3-connected, $|L(v_2)|\geq 3$. Choose two colors $a$ and $b$ from $L(v_2)$. Delete $a$ and $b$ from the lists of all neighbors of $v_2$ in $B_7$ and let $L'(v_2)=\{a,b\}$. After the coloring of $G[B_7]$ we have $2=|L'(v_2)|>d_H(v_2)=1$, so we are done by (ii).

This completes the proof of the theorem.
\hfill $\Box$ \msk

 \noindent
{\bf Remark:} The above proof also works for $k=6$, except in the
final case.

\section{The case $k=5$}

\begin{theo}
Let\/ $k=5$. There is a planar, 3-connected, non-complete
graph\/ $G$ with\/ $d(S_5)=4$ and a list assignment\/ $L$ with\/
$|L(v)|=\min \{d(v), 5\}$ for all\/ $v\in V$ such that\/ $G$ is
not\/ $L$-list colorable.
\end{theo}


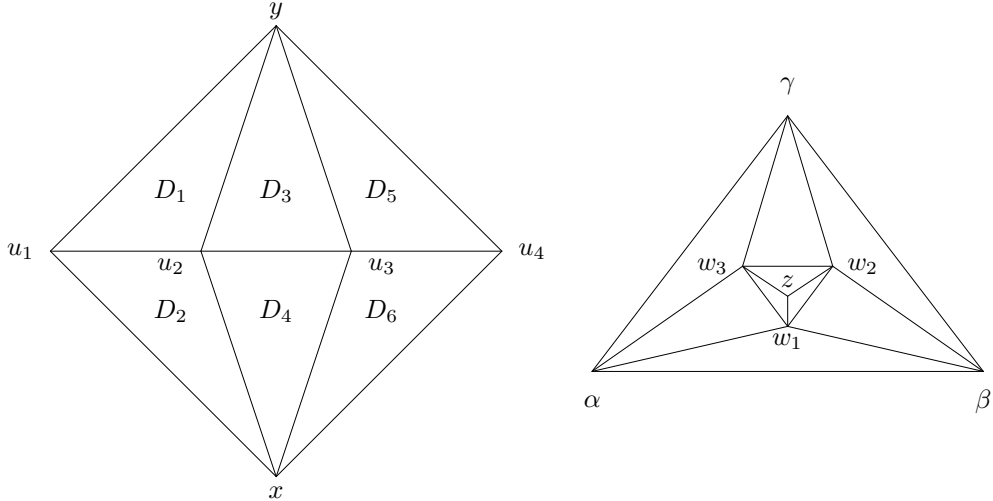
\begin{figure}[h]
\unitlength1mm
\hspace{-5mm}
\begin{picture}(148,76)
\thinlines
\drawpath{42.0}{68.0}{12.0}{38.0}
\drawpath{12.0}{38.0}{42.0}{8.0}
\drawpath{42.0}{8.0}{72.0}{38.0}
\drawpath{72.0}{38.0}{42.0}{68.0}
\drawpath{42.0}{68.0}{32.0}{38.0}
\drawpath{32.0}{38.0}{42.0}{8.0}
\drawpath{42.0}{8.0}{52.0}{38.0}
\drawpath{52.0}{38.0}{42.0}{68.0}
\drawpath{12.0}{38.0}{32.0}{38.0}
\drawpath{32.0}{38.0}{52.0}{38.0}
\drawpath{52.0}{38.0}{72.0}{38.0}
\drawcenteredtext{42.0}{6.0}{$x$}
\drawcenteredtext{42.0}{70.0}{$y$}
\thicklines
\drawcenteredtext{8.0}{38.0}{$u_1$}
\drawcenteredtext{28.0}{36.0}{$u_2$}
\drawcenteredtext{56.0}{36.0}{$u_3$}
\drawcenteredtext{76.0}{38.0}{$u_4$}
\drawcenteredtext{28.0}{46.0}{$D_1$}
\thinlines
\drawcenteredtext{28.0}{30.0}{$D_2$}
\drawcenteredtext{42.0}{46.0}{$D_3$}
\drawcenteredtext{42.0}{30.0}{$D_4$}
\drawcenteredtext{56.0}{46.0}{$D_5$}
\drawcenteredtext{56.0}{30.0}{$D_6$}
\drawpath{84.0}{22.0}{110.0}{56.0}
\drawpath{116.0}{36.0}{104.0}{36.0}
\drawpath{136.0}{22.0}{84.0}{22.0}
\drawpath{84.0}{22.0}{104.0}{36.0}
\drawpath{104.0}{36.0}{110.0}{56.0}
\drawpath{110.0}{56.0}{116.0}{36.0}
\drawpath{116.0}{36.0}{136.0}{22.0}
\drawpath{136.0}{22.0}{110.0}{28.0}
\drawpath{110.0}{28.0}{84.0}{22.0}
\drawpath{104.0}{36.0}{110.0}{28.0}
\drawpath{110.0}{28.0}{116.0}{36.0}
\drawcenteredtext{84.0}{18.0}{$\alpha$}
\drawcenteredtext{110.0}{26.0}{$w_1$}
\drawcenteredtext{120.0}{36.0}{$w_2$}
\drawcenteredtext{100.0}{36.0}{$w_3$}
\drawpath{110.0}{56.0}{136.0}{22.0}
\drawcenteredtext{136.0}{18.0}{$\beta$}
\drawcenteredtext{110.0}{60.0}{$\gamma$}
\drawpath{104.0}{36.0}{110.0}{32.0}
\drawpath{110.0}{32.0}{116.0}{36.0}
\drawpath{110.0}{32.0}{110.0}{28.0}
\drawcenteredtext{110.0}{34.0}{$z$}
\end{picture}
\vspace{-1cm}
\caption{(a) Subgraph $H$.~~(b) The structure inside each $D_i$.}
\end{figure}

{\bf Proof.} Consider the subgraph $H$ in Figure 2a with list
assignment $L'(x)=\{a\}$, $L'(y)=\{b\}$, and $L'(u_i) =
\{a,b,1,2,3\}$, for $i=1,\dots,4$. We claim that  three colors
can be assigned to each of the triangles $D_1,\dots, D_6$ such
that in every coloring of $H$ at least one of the triangles is
colored by its assigned colors; note that the order of the colors
in each assignment matters.  The assignments are
$$D_1: 2,1,a; \ \ D_2: 3,1,b; \ \ D_3: 2,3,a; \ \ D_4: 3,2,b; \ \ D_5:
1,2,a; \ \ D_6: 1,3,b.$$ Note that the prescribed coloring of
$D_1$ or $D_2$ occurs if $u_2$ is colored by 1. Analogously, the
prescribed coloring of $D_5$ or $D_6$ occurs if $u_3$ is colored
by 1. Thus we may assume that $u_2$ and $u_3$ are not colored by
1. But now the prescribed coloring of $D_3$ or $D_4$ occurs.

Now assume that inside of each of $D_1,\dots,D_6$ we have the
structure in Figure 2b. For each triangle $D_i$, let
$\alpha,\beta,\gamma$ be the colors assigned to $D_i$. Let
$L'(w_1)=\{\alpha,\beta,4,5,6\}, L'(w_2)=\{\beta,\gamma,4,5,6\}$
and $L'(w_3)=\{\alpha,\gamma,4,5,6\}$ and $L'(z)=\{4,5,6\}$. Thus
$H$ is not $L'$-list colorable for the given list assignment.

Now build a graph $G$ by taking 25 copies of $H$ and identifying
all 25 $x$-vertices to a vertex $x^*$ and all 25 $y$-vertices to
vertex $y^*$. Join vertex $u_4$ in the $i$-th copy with vertex
$u_1$ in the $(i+1)$-th copy, for $i=1,\dots, 24$. Define a list
assignment $L$, as follows. Let $L(x^*)=\{7,8,9,10,11\}$ and
$L(y^*)=\{12,13,14,15,16\}$. For all other vertices let
$L(v)=L'(v)$, where each pair of $ L(x^*)\times L(y^*)$
corresponds with the color pair $(a,b)$ of one copy of $H$.

Finally, observe that $G$ is a planar, 3-connected non-complete
graph with a list assignment $L$ such that $|L(v)|=\min\{d(v),5\}$
for all $v\in V$, but $G$ is not $L$-list colorable. \hfill $\Box$

\msk

Note that $d(S_5)=4$ for the example in the proof above. In fact,
it is almost the same example as in \cite{TuzaV02}.

\section{Open problems}

Despite our progress in this paper, two essential problems remain open.

\begin{prob}
Let\/ $k\geq 6$ be an integer and let\/ $G$ be a non-complete
planar graph with\/ $\kappa(G)\in\{3,4\}$ and\/ $d(S_k)=2$. Is\/
$G$ $f$-list colorable when\/ $f(v)=\min\{d(v),k\}$ for all\/
$v\in V$? \ssk
\end{prob}

\begin{prob}
Let\/ $k=5$ be an integer and let\/ $G$ be a connected planar
graph that is not a Gallai tree and that has\/ $d(S_k)\geq 5$.
Is\/ $G$ $f$-list colorable when\/ $f(v)=\min\{d(v),k\}$ for all\/
$v\in V$?
\end{prob}

One can raise the analogous questions for $K_5$-minor-free graphs,
too. We know from Proposition~\ref{c3,d2} that the answer to the
first one with $\kappa=3$ is negative.


\def \ncm  {\newcommand}

\def \bibi {\bibitem}
\def \etal {{\em et al.}}

\ncm{\au}[1]{{\sc #1\,:\,}}      \ncm{\ti}[1]{{\it #1\/.\,}}
\ncm{\bkti}[1]{{\sl #1}}         \ncm{\jou}[1]{{\rm #1\/}}
\ncm{\vol}[1]{{\bf #1}}          \ncm{\yr}[1]{\rm (#1),}
\ncm{\pp}[1]{{\rm #1.}}

\def \JCTB {Journal of Combinatorial Theory, Ser.~B}
\def \DM {Discrete Mathematics}
\def \DMGT {Discussiones Mathematicae Graph Theory}
\def \JGT {Journal of Graph Theory}
\def \CoNum {Congressus Numerantium}

\end{document}